\documentclass[12pt,psamsfonts]{amsart}
\usepackage{amsmath}
\usepackage{amsthm}
\usepackage{amssymb}
\usepackage{amscd}
\usepackage{amsfonts}
\usepackage{amsbsy}
\usepackage[latin1]{inputenc}
\usepackage{epsfig,afterpage}
\usepackage[dvips]{psfrag}
\usepackage{psfrag}
\usepackage[all]{xy}
\usepackage{color}
\usepackage{multirow}
\usepackage{pdflscape}
\usepackage{slashbox} 
\usepackage{pict2e} 
 \usepackage{array}

\newcommand{\R}{\ensuremath{\mathbb{R}}}

\newcommand{\e}{\varepsilon}

\newcommand{\sgn}{\mbox{\normalfont sgn}}

\usepackage{overpic}

\newtheorem {theorem} {Theorem}
\newtheorem* {proposition}  {Fundamental Lemma}

\newtheorem {remark} {Remark}

\AtBeginDocument{%
  
}

\begin{document}

\title[Invariant manifolds for piecewise smooth systems]{Invariant manifolds for piecewise smooth differential systems in $\mathbb{R}^{3}$}

\author[Buzzi, C.A., Euz\'{e}bio, R.D. and Mereu, A.C.]
{Claudio A. Buzzi$^{1}$, Rodrigo D. Euz\'{e}bio$^{2}$ and Ana C. Mereu$^3$}

\address{$^1$ UNESP--IBILCE, CEP 15054--000, S. J. Rio Preto, S\~ao Paulo, Brazil}
\email{claudio.buzzi@unesp.br}

\address{$^2$ Departamento de Matem\'{a}tica, IME--UFG, R. Jacarand\'{a}, Campus Samambaia, Zip Code 74001-970, Goi\^{a}nia, GO, Brazil}
\email{euzebio@ufg.br}

\address{$^3$ Departamento de F\'{i}­sica, Qu\'{i}mica e Matem\'{a}tica.
 UFSCar. 18052-780, Sorocaba, SP, Brazil}
\email{anamereu@ufscar.br}



\date{}
\dedicatory{}

\begin{abstract}
In this paper some piecewise smooth perturbations of a three-dimensional differential system are considered. The existence of invariant manifolds filled by periodic orbits is obtained after suitable small perturbations of the original differential system. These manifolds emerge from a continuum of cylinders of $\mathbb{R}^3$ which does exist for the piecewise smooth differential systems after a rotation of some planar algebraic polynomial curves. The main tool used in order to obtain the results is the averaging theory for piecewise smooth differential systems.
\end{abstract}

\maketitle


\section*{Introduction}\label{secao introducao}

The existence of invariant manifolds plays an important role in the qualitative understanding of trajectories in a dynamical systems. They occur, for instance, in Hamiltonian system presenting a set of variables known as {\it action-angle variables}. In these variables the equations of motion are very simple and when the trajectories are bounded, a precession of the angle variables around a torus takes place. We also address the KAM theory (see \cite{Arn1963,Kol1954,Mos1962}) which basically consists in answering the question: ``\textit{Do the trajectories of the perturbed Hamiltonian system still lie on invariant tori, at least for small perturbations?}''. Accordingly, to detect how persistent are some objects is actually one of the main problems of general perturbation theory. This problem is particularly interesting when the perturbation involves the persistence of invariant manifolds because the dynamics is confined to a lower dimension object so usually the study is simpler.

In this paper we consider perturbations of a three-dimensional differential system having a continuum of invariant cylinders. We are mainly motivated by the problem of finding limit cycles in dimension two, from a perturbative point of view. Indeed, in dimension two, it is very effective to perturbate a linear center in order to obtain limit cycles, so in our context we consider a three-dimensional differential system behaving like a linear center but the extra dimension generates a continuum of invariant cylinders for the model we are considering. A simplified differential system achieving these assumptions can be taken as the following
\begin{equation}
\begin{array}{l}\label{system_no_perturbation}
\dot{x}=-y,\vspace{0.2cm}\\
\dot{y}=x,\vspace{0.2cm}\\
\dot{z}=h(x,y).
\end{array}
\end{equation}
where $h$ is a continuous real-valued function. The main goal of this paper is to study the existence of two dimensional invariant manifolds under small piecewise smooth perturbations of 
\eqref{system_no_perturbation}.

Invariant manifolds in smooth and piecewise smooth systems have been recently studied in some papers. Concerning smooth systems, the existence of invariant manifold were addressed in \cite{LliRebTor2011a,LliRebTor2011b,LliTei2009} among others. In the piecewise linear context, several authors study the existence, the number and the stability of invariant cones and limit cycles in piecewise linear systems, see \cite{CarmFernFrei2011,CarFrePonTor2005,huan2017,huan2016} and references therein. For the case of nonlinear piecewise smooth systems similar results are obtained in \cite{KuppHosh2011,WeisKuppHosh2012,WeisKuppHosh2015} where the considered systems are approximated by linear ones. In this case, the presence of the so called {\it sliding motion} have also been addressed.

We stress that the study of general nonlinear piecewise smooth differential systems can be a hard task. This is the main goal of this paper. Regardless, some particular situations have been considered in the literature. For instance, in \cite{BuzEuzMer2016} the authors study smooth and piecewise smooth perturbations of a cylinder filled by periodic orbits obtaining the persistence of limit cycles. A similar approach is done in \cite{CarvFrei2017} where the authors study the existence of limit cycles around a fold-fold singularity in $\mathbb{R}^3$. The study of invariant manifolds have been also done in \cite{LlibMartTeix2010,LlibMereTeix2010}. In these papers the authors provide sufficient conditions for the existence of some particular invariant manifolds, namely cylinders and tori.

This paper is organized as follows: Section \ref{secao colocacao do problema} is devoted to introduce the problem, notations and methods used throughout the paper. In Section \ref{results} we present the main results of the paper and their proofs are provided in Section \ref{provas_discussoes}. In Section \ref{Examples} some particular examples are discussed and finally Section
\ref{section_averaging} provides the averaging method for piecewise smooth differential systems. 

From now on we refer to the differential system \eqref{system_no_perturbation} simply by system \eqref{system_no_perturbation} dropping the word differential.

\section{Setting the problem}\label{secao colocacao do problema}

We start noticing that solutions of system \eqref{system_no_perturbation} starting at $(x_0,y_0,z_0)\in\mathbb{R}^3$ are $x(t)= x_0\cos t-y_0\sin t$, $y(t)= y_0\cos t+x_0\sin t$ and

\begin{equation}\label{solution_z(t)}
z(t)= z_0+\displaystyle\int_{0}^{t}h(x_0\cos t-y_0\sin t,y_0\cos t+x_0\sin t)dt.
\end{equation}

We are interested in the study of the cylinders $\textsl{C}_{\rho}=\{(x,y,z)\in \R^3:x^2+y^2 =\rho^2\}$, which are invariant sets for system \eqref{system_no_perturbation} for all values $\rho\neq0$ regardless of $h$. To see that notice that $x(t)^2+y(t)^2=\rho_0^2$ for any initial condition satisfying $x_0^2+y_0^2=\rho_0^2$ with $\rho_0\neq0$. Moreover, the invariance of the cylinders $\textsl{C}_{\rho}$ does not depend on the function $h$. When $\rho_0=0$ we have initial conditions on the $z-$axis which is an invariant set for the flow associated to system \eqref{system_no_perturbation}. For our purposes this case is not considered.

In this paper we shall assume that the cylinders $\textsl{C}_{\rho}$ are periodic, that is, every trajectory contained on it is periodic having the same period. Of course, in order to do that we must assume extra hypothesis on the function $h$.

\begin{remark}
	The periodicity of the cylinders is mandatory for obtaining the results of this paper. For instance, if $h(x,y)=k$, where $k\in\mathbb{R}\setminus\{0\}$, then cylinders $\textsl{C}_{\rho}$ are still invariant and trajectories spiral around them and go to infinity. However, the methods used in this paper do not apply in such case, see Section \ref{section_averaging}.
\end{remark}

In order to obtain a periodic behavior on the cylinder is sufficient to choose the function $h$ properly such that: $(i)$ the solution $z(t)$ presented in \eqref{solution_z(t)} is periodic and $(ii)$ the averaging method applies. This is achieved assuming, respectively, the following hypotheses:

\begin{itemize}
	\item[{\bf H1:}] $h(x,y)=\varphi(x^2+y^2)\overline{h}(x,y)\label{condicao_funcao_h}$,
	 $\overline{h}(x,y)=x\,P_1(x^{2},y^{2})+xy\,P_2(x^{2},y^{2})+y\,P_3(x^{2},y^{2})$,
with $P_1$, $P_2$ and $P_3$ polynomial functions and $\varphi$ is an arbitrary continuous function.
\end{itemize}

\begin{itemize}
	\item [{\bf H2:}] The function $h$ expressed in terms of cylindrical coordinates does not depend on the radius $r$, that is, $\dfrac{\partial h}{\partial r}(r,\theta)\equiv 0$.
\end{itemize}

Hypothesis ${\bf H1}$ is already considered in \cite{BuzEuzMer2016}. Now, under the previous hypotheses, to study the appearance of periodic invariant manifolds from system \eqref{system_no_perturbation}, we will perform a small \textit{radial cylindrical perturbation} of this system. In the literature (see \cite{LliRab2015}) a center is rigid if its angular speed is constant. Here we say that a vector field on $\R^3$ is radial cylindrical if in cylindrical coordinates $x=r\cos\theta,$ $y=r\sin\theta$ and $z=z$ we have $\dot{\theta}=0$ and $\dot{z}=0$. More specifically, this perturbation is given by the function
$$
g(x,y,z)=
\left\{
\begin{array}{rl}
g^+(x,y,z),&\quad y>0,\vspace{0.3cm}\\
g^-(x,y,z),&\quad y<0,
\end{array}
\right.
$$
where $g^{\pm}=(g_1^{\pm},g_2^{\pm},g_3^{\pm})$ are given by
\begin{equation}\label{perturbations}
\begin{array}{rcl}
g_1^{\pm}(x,y,z)&=&x\Psi^\pm(x,y,z),\\
g_2^{\pm}(x,y,z)&=&y\Psi^\pm(x,y,z),\\
g_3^{\pm}(x,y,z)&=&0,
\end{array}
\end{equation}
and
\[\Psi^\pm(x,y,z)=\displaystyle\sum_{0\leq i+j+k\leq n-1} a^{\pm}_{ijk}x^{i}y^{j}z^{k},\]
with $i,j,k,n\in\mathbb{N}$ and $a_{ijk}\in\mathbb{R}$, $\forall i,j,k\in\mathbb{N}$. Assuming  $X=(x,y,z)$, that lead us to the piecewise smooth differential system
\begin{equation}\label{system_noncontinuous}
\dot{X}=f(X)+\varepsilon g(X),
\end{equation}
where $\varepsilon$ is a small parameter and $f$ is given by $f(X)=(-y,x,h(x,y))^T$.

System \eqref{system_noncontinuous} with $\varepsilon=0$ is the smooth system presenting the con\-ti\-nuum of cylinders mentioned previously, but assuming $\varepsilon\neq 0$ sufficiently small one cannot assure they are preserved. Moreover, the plane $y=0$ split the original cylinders into two parts, so their persistence depend also on each perturbation is taking place (that is, if $g^+$ or $g^-$ is acting).

Notice that over the plane $\Sigma=\{(x,0,z)\in\R^3\}$ the differential system \eqref{system_noncontinuous} is bi-valued therefore trajectories on this plane must be properly defined. According to \cite{Fil1988} such trajectories are of {\it crossing}, {\it sliding} or {\it tangent} type depending on how the vector fields $f(X)+\varepsilon g^{\pm}(X)$ interacts to $\Sigma$. The first (respect. second) case occurs when for points on $\Sigma$ the vectors fields point to the same side of it (respect. point to opposite directions). The third case occurs when at least one the vector fields is tangent to $\Sigma$. For the particular case addressed in the paper, it is easy to see that every point on $\Sigma\setminus\{x=0\}$ is of crossing type. The set $\Sigma\cap\{x=0\}$ (the $z-$axis) is formed by tangency points of both vector fields. For more details on general properties of piecewise smooth differential systems see \cite{BerBudChaKow2008,Fil1988}.

Here in this paper we are interested in periodic orbits of crossing type. Under this assumption and assuming hypothesis {\bf H1} and {\bf H2}, we can apply the averaging theory for piecewise smooth systems in order to obtain the bifurcation of two dimensional manifolds desired. See Theorem \ref{averaging} of Section \ref{section_averaging} which can be found originally in \cite{LliNovTei2015c}.

\section{Statements of the main results}\label{results}

Now we establish the main results of the paper. For that, we call
$$
I_h(\theta)=\displaystyle\int_0^\theta h(\cos s,\sin s)ds.
$$

Define also the integrals

$$
c^0_{10}=\int_0^\pi I_\theta \,d\theta,\quad
c^0_{01}=\int_\pi^{2\pi} I_\theta \,d\theta,\quad
c^1_{11}=\int_0^\pi I_\theta C_\theta \,d\theta,\quad
$$

$$
c^1_{12}=\int_\pi^{2\pi} I_\theta C_\theta \,d\theta, \quad
c^1_{21}=\int_0^\pi I_\theta S_\theta \,d\theta, \quad
 c^1_{22}=\int_\pi^{2\pi} I_\theta S_\theta \,d\theta,\quad
$$

$$
c^2_{10}=\int_0^\pi I_\theta^2 \,d\theta,\qquad
c^2_{01}=\int_\pi^{2\pi} I_\theta^2 \,d\theta,
$$
where $I_\theta=I_h(\theta)$, $C_\theta^i=\cos^i(\theta)$ and $S^i_\theta=\sin^i{\theta}$, $i=1,2$. A key result of the paper follows as an application of averaging theorem.

\begin{proposition}\label{t1}
Consider system \eqref{system_noncontinuous} and assume that hypotheses {\bf H1} and {\bf H2} hold for some function $h$ and $g$ is radial cylindrical. Consider also cylindrical coordinates $x=r\cos\theta$, $y=r\sin\theta,$ $z=z$, with $r>0$ and suppose that the real valued function
$$
\begin{array}{rcrl}
\psi^n_{z_0}(r)&=&r\left[\displaystyle\int_{0}^{\pi} \Psi^{+}(r\cos\theta,r\sin\theta,z_0+I_h(\theta))\; d\theta \right.\vspace{0.3cm}\\
&+&\left.\displaystyle\int_{\pi}^{2\pi} \Psi^{-}(r\cos\theta,r\sin\theta,z_0+I_h(\theta))\; d\theta \right]
\end{array}
$$
has a zero $r_0=r_0(z_{0})$ such that $(\partial \psi^n_{z_0}/\partial r)(r_0(z_0))\neq0$. Call $S=\{z_0\in \mathbb{R};(\partial \psi^n_{z_0}/\partial r)(r_0(z_0))\neq0\}$ and define the function $\psi_n:\mathbb{R}_+\times S\longrightarrow\mathbb{R}$ by $\psi_n(r,z_0)=\psi^n_{z_0}(r)$. Then, for $|\e|$ sufficiently small and for each fixed $z_0\in S$, the following statements hold:
\begin{itemize}
	\item[(i)] there exists a $2\pi$--periodic solution $X(t,z_{0},\e)$ of system \eqref{system_noncontinuous} such that $X(0,z_{0},\e)\to (r_0(z_0),0,z_0)$ as $\e\to 0$.
	\item[(ii)] the trajectory of each point starting on $\psi_n^{-1}(0)$ is $2\pi$-- periodic. In other words, system \eqref{system_noncontinuous} has a piecewise smooth manifold $\mathfrak{M}$ homeomorphic to the revolution manifold $\widetilde{\mathfrak{M}}=\mathbb{S}^{1}\times \psi_n^{-1}(0)\subset\mathbb{R}^3$ filled by closed orbits.
\end{itemize}
\end{proposition}

\begin{remark}
	We notice that we only consider the derivative of function $\psi^n_{z_0}$ in Fundamental Lemma instead of deal with the Brouwer degree, as considered in Theorem \ref{averaging}. The reason why we do that is because in the most part of problems with a polynomial perturbation approach, the averaged function $\psi^n_{z_0}$ is also polynomial. In this case, the study of Brouwer's degree is equivalent to the study of the signal of the derivative of $\psi^n_{z_0}$.
\end{remark}

We stress that Fundamental Lemma is quite general, since the degrees of the perturbation $g$ of system \eqref{system_noncontinuous} are arbitrary and so is function $h$. However, it does not provide any information about the shape of the manifolds bifurcating from the continuum of cylinders. In this direction, the next results are particular cases of Fundamental Lemma which provide the nature of the manifolds $\mathfrak{M}$ depending on the degrees of the perturbations. It also assure the realization of such manifolds.

\begin{theorem}\label{teo_generalv}
Assume that system \eqref{system_noncontinuous} satisfies hypotheses ${\bf H1}$ and ${\bf H2}$ for some function $h$ and $g$ is a radial cylindrical perturbation of degree $n$. Then function $\psi_n:\mathbb{R}_+\times S\rightarrow\R$ writes 
	\begin{equation}\label{poly}
	\psi_n(r,z_0)=r\left[\displaystyle\sum_{0\leq i+j\leq n-1}C^n_{ij}r^iz_0^j\right],
	\end{equation} 
	where the coefficients $C^n_{ij}$ are real constants depending on both degree and parameters of the perturbation $g$. Consequently we get:
	\begin{itemize}
		\item[(i)]  $\psi_n^{-1}(0)$ is the intersection of an algebraic real curve of degree $n-1$ in the $(r,z_0)$-plane with the set $\mathbb{R}_+\times S$;
		\item[(ii)] $\mathfrak{M}$ is homeomorphic to the revolution manifold $\widetilde{\mathfrak{M}}$ obtained by the rotation of $\psi_n^{-1}(0)$ around the $z_0-$axis.
	\end{itemize}
 Moreover, given any polynomial $\bar{\psi}_n(r,z_0)$, as in \eqref{poly} of degree $n$, there exists a perturbation $g$, also of degree $n$, that realizes the polynomial $\bar{\psi}_n$.
\end{theorem}

The following results are particular cases of previous theorem when we consider perturbations of degree one and two, respectively.

\begin{theorem}\label{teo1}
Under the same hypotheses of Theorem \ref{teo_generalv}, with $n=2$, we have $\psi_2(r,z) =r\left[ C^2_{10} r + C^2_{01} z + C^2_{00}\right]$, where
\[
\begin{array}{l}
C^2_{10}=2(a_{010}^{+}-a_{010}^{-}),\vspace{0.2cm}\\
C^2_{01}=\pi(a_{001}^{+}+a_{001}^{-}),\vspace{0.2cm}\\
C^2_{00}=\pi( a_{000}^{+}+a_{000}^{-})+c^0_{10} a_{001}^{+}+c^0_{01}a_{001}^{-}.
\end{array}
\]
Moreover, if we call $L_1$ the straight line $\{(r,z):\psi_2(r,z) =0\mbox{ and } r>0\}$ in the $(r,z)$-plane and $\ell=L_1\cap\{\theta=0\}$ then we have
\begin{itemize}
	\item[(i)]  if $C^2_{10}\neq0$ then  $\widetilde{\mathfrak{M}}$ is a revolution manifold surrounding the $z-$axis obtained by revolution of the segment $\ell$.
	\item[(ii)] if $C^2_{10}=0$ then the method does not apply and we can not obtain manifolds fulfilled by periodic orbits. 
\end{itemize}
\end{theorem}

\begin{theorem}\label{teo2}
		Under the same hypotheses of Theorem \ref{teo_generalv}, with $n=3$, we have $\psi_3(r,z) = r\left[C^3_{20} r^2 + C^3_{11} rz +C^3_{02} z^2 +C^3_{10} r +C^3_{01} z + C^3_{00}\right]$, where
		$$
		\begin{array}{rcl}
		C^3_{20}&=&\dfrac{\pi}{2}(a^+_{200}+a^+_{020}+a^-_{200}+a^-_{020}),\vspace{0.2cm}\\
		C^3_{11}&=&2(a^+_{011}-a^-_{011}),\vspace{0.2cm}\\
		C^3_{02}&=&\pi(a^+_{002}+a^-_{002}),\vspace{0.2cm}\\
		C^3_{10}&=&a^+_{101}c^1_{11}+a^+_{011}c^1_{21}+2(a^+_{010}-a^-_{010})+a^-_{101}c^1_{12}+a^-_{011}c^1_{22},\vspace{0.2cm}\\
		C^3_{01}&=&2a^+_{002}c^0_{10}+\pi a^+_{001}+2a^-_{002}c^0_{01}+\pi a^-_{001},\vspace{0.2cm}\\
		C^3_{00}&=&a^+_{002}c^2_{10}+a^+_{001}c^0_{10}+\pi(a^+_{000}+a^-_{000})+a^-_{002}c^2_{01}+a^-_{001}c^0_{01}.
		\end{array}
		$$
Moreover, $\psi_3^{-1}(0)\setminus\{r=0\}$ is the intersection of a conic with $\mathbb{R}_+\times S$. So, when $\psi_3^{-1}(0)\setminus\{r=0\}\neq\emptyset$, $\widetilde{\mathfrak{M}}$ is a manifold of revolution obtained by the rotation of (eventually only a piece of) a conic around the $z-$axis. 
\end{theorem}

We stress that while in Theorem \ref{teo1} we describe the manifolds $\widetilde{\mathfrak{M}}$, in Theorem \ref{teo2} they are a large number of objects in the sense that a full classification is an arduous topological task and is out of the goal of this paper. For instance, it is not difficult to see that such objects must include manifolds with more than one connected component not always compact (in the case, for instance, of a hyperbola). 

\section{Proofs of the main results}\label{provas_discussoes}

Now we apply the methods and tools described previously in order to prove the results presented in Section \ref{results}.

\subsection{Proof of Fundamental Lemma}

Consider system \eqref{system_noncontinuous} and assume hypothesis ${\bf H1}$, so solutions of such system with $\varepsilon=0$ are periodic. Moreover they live on the cylinders $\textsl{C}_{\rho}$, then we perform a cylindrical change of coordinates in system \eqref{system_noncontinuous} by introducing the new variables $(z,r,\theta)$ given implicitly by $x=r\cos\theta$, $y=r\sin\theta$ and $z=z$. In the new variables $(z,r,\theta)$ system \eqref{system_noncontinuous} writes
	\begin{equation}\label{system_cylindrical_coordinates}
	\begin{array}{rl}
	\dot{r}=&\varepsilon r \Psi^{\pm}(r\cos\theta,r\sin\theta,z) ,\vspace{0.2cm}\\
	\dot{\theta}=&1,\vspace{0.2cm}\\
	\dot{z}=&h(r\cos\theta,r\sin\theta).\vspace{0.2cm}
	\end{array}
	\end{equation}
Observe that, in cylindrical coordinates, the discontinuity plane is provided by $\sgn(\sin\theta)$. Now we change the independent variable $t$ of system \eqref{system_cylindrical_coordinates} to the new variable $\theta$. By the hypothesis ${\bf H2}$ it follows that $h(r\cos\theta,r\sin\theta)=h(\cos\theta,\sin\theta)$ so we obtain the solution of third equation of \eqref{system_cylindrical_coordinates}, namely, $z(\theta)=z_{0}+I_h(\theta)$. Consequently, replacing this expression into the first equation of system \eqref{system_cylindrical_coordinates} we get the equation
\begin{equation}
\begin{array}{rl}\label{final_system}
\dfrac{dr}{d\theta}=\varepsilon r \Psi^{\pm}\left(r\cos\theta,r\sin\theta,z_{0}+I_h(\theta)\right).
\end{array}
\end{equation}
Finally, in order to apply Theorem \ref{averaging}, we observe that system \eqref{final_system}
can be written as
\begin{equation}
\begin{array}{rl}\label{final_system2}
\dfrac{dr}{d\theta}=\varepsilon G(\theta,r) +\varepsilon^{2}R(\theta,r,\e),
\end{array}
\end{equation}
with
\[
\begin{aligned}
&G(\theta,r)=G^+(\theta,r)+\sgn(\eta(\theta,r))G^-(\theta,r),\\
&R(\theta,r,\e)=R^+(\theta,r,\e)+\sgn(\eta(\theta,r))R^-(\theta,r,\e),
\end{aligned}
\]
where
\[
\begin{aligned}
&\eta(\theta,r)=\sin \theta,\\
&G^+(\theta,r)=(\Psi^+(\vartheta)+\Psi^-(\vartheta))/2,\\
&G^-(\theta,r)=(\Psi^+(\vartheta)-\Psi^-(\vartheta))/2,\\
&R^+(\theta,r,\e)=0,\\
&R^-(\theta,r,\e)=0,\\
&\vartheta=(r\cos\theta,r\sin\theta,z_0+I_h(\theta)).
\end{aligned}
\]
It is clear that the functions $G^\pm:\R\times\mathbb{R}_+\rightarrow\R$, $R^\pm:\R\times\mathbb{R}_+\times(-\e_0,\e_0)\rightarrow\R$, $\eta:\R\times\mathbb{R}_+\rightarrow\R$ are continuous, $2\pi$-periodic in the variable $\theta$, locally Lipschitz with respect to $r$ and $\eta$ is a $C^1$ function having $0$ as a regular value. From Remark \ref{rem1} the condition $(ii)$ of Theorem \ref{averaging} holds, because $\partial_\theta\eta(0,r)=\cos(0)=1$ and $\partial_\theta\eta(\pi,r)=\cos(\pi)=-1$. Now we call $\psi^n_{z_0}:\mathbb{R}_+\rightarrow\R$ given by
\begin{equation}\label{function_F}
\psi^n_{z_0}(r)= \displaystyle \int_{0}^{2\pi}G(\theta,r)d\theta,
\end{equation}\label{psi_fundamental_lemma}
where $n$ denotes the degree of the perturbation function $g$. Therefore function $\psi^n_{z_0}$ is polynomial in the variable $r$ and writes
\begin{equation*}
\begin{array}{rl}
\psi^n_{z_0}(r)=&r\left[\displaystyle \int_{0}^{\pi}\Psi^+\left(\vartheta\right)d\theta+\displaystyle \int_{\pi}^{2\pi}\Psi^-\left(\vartheta\right)d\theta\right].
\end{array}
\end{equation*}

Therefore, by hypothesis, $\psi^n_{z_0}$ has a zero $r_0=r_0(z_{0})$ satisfying \linebreak $(\partial \psi^n_{z_0}/\partial r)(r_0(z_0))\neq0$. Consequently, $z_0\in S$ and then there exist a neighborhood $V_{r_0} \subset\mathbb{R}_+$ of $r_0(z_0)$ such that $\psi^n_{z_0}(r)\neq 0$ for all $r\in\overline{V_{r_0}}\backslash\{r_0(z_0)\}$. Moreover, by the late, $$
d_B(\psi^n_{z_0},V_{r_0},r_0(z_0))=\sgn(\partial \psi^n_{z_0}/\partial r)(r_0(z_0)),
$$
that is $d_B(\psi^n_{z_0},V_{r_0},r_0(z_0))$ is $+1$ or $-1$. Thus bullet $(iii)$ of Theorem \ref{averaging} holds. 

Then, by applying Theorem \ref{averaging} for the $z_{0}$-parametric system \eqref{final_system}, it follows that, for $|\varepsilon|>0$ sufficiently small, system \eqref{final_system} has a $2\pi$-periodic orbit $r(\cdot,z_0,\varepsilon)$ for all $z_{0}\in S$ such that $r(0,z_0,\varepsilon)>0$ and such that $r(0,z_0,\varepsilon)\to r_0(z_0)$ when $\varepsilon\to 0$. The same is true for system \eqref{system_cylindrical_coordinates} for every $z_0\in S$, since they are equivalent to system \eqref{final_system}. Therefore, coming back through the cylindrical coordinates $x=r\cos\theta$, $y=r\sin\theta$, $z=z$,  there exists a $2\pi$--periodic solution $X(t,z_{0},\e)$ of system \eqref{system_noncontinuous} such that $X(0,z_{0},\e)\to (r_0(z_0),0,z_0)$ as $\e\to 0$, thus the proof of the first statement of the Fundamental Lemma is done. 

In order to prove the second statement of Fundamental Lemma, observe that since $\psi_n$ is defined only for $z_0\in S$, then for each pair $(r_0(z_0),z_0)$ we have $\psi_n(r_0,z_0)=\psi^n_{z_0}(r_0)=0$ so we can apply bullet $(i)$ proved before to guarantee that such a pair corresponds to a $2\pi$-periodic solution. That is to say we have $2\pi$-periodic solutions for every initial condition on $\psi_n^{-1}(0)$. Again, coming back from cylindrical to Cartesian coordinates, every point on the curve $\psi_n(r_0,z_0)=0$, whose trace is contained on $\mathbb{R}_+\times S$, corresponds to a closed orbit $\Gamma_{z_0}$ passing through $(r(z_0),0,z_0)$, being $\Gamma_{z_0}$ contained in $\mathbb{R}^3$. More specifically, such an orbit is formed by two pieces $\Gamma_{z_0}^{\pm}$ corresponding to each crossing through the plane $y=0$: $\Gamma_{z_0}^+=[0,\pi]\times(r_0,z_0)$ and $ \Gamma_{z_0}^-=[\pi,2\pi]\times(r_0,z_0).$ Therefore, $\Gamma_{z_0}=\Gamma_{z_0}^{+}\bigcup\Gamma_{z_0}^{-}$. Proceeding analogously for every $z_0\in S$ we get the revolution manifold $\widetilde{\mathfrak{M}}=\mathbb{S}^{1}\times \psi_n^{-1}(0)$. It is easy to see then that
$$
\mathfrak{M}=[0,\pi]\times \psi_n^{-1}(0)\bigcup\;[\pi,2\pi]\times \psi_n^{-1}(0)
$$
is homeomorphic to $\widetilde{\mathfrak{M}}$ due to the eventually non regular contact with the plane $y=0$. Moreover, $\mathfrak{M}$ is filled by closed orbits of system \eqref{system_noncontinuous}. This concludes the proof of Fundamental Lemma.

\subsection{Proof of Theorem \ref{teo_generalv}}

For the first part of the proof we consider $\Psi^\pm$ a polynomial of degree $n-1$. So, the functions $g_i^{\pm}(\vartheta)$, $i=1,2,3$, in \eqref{system_cylindrical_coordinates} are polynomial of degree $n$ in the variables $r,z$. It implies that the function  $\psi_{z_0}^n$, defined in \eqref{function_F} also is polynomial of degree $n$ in the variables $r$ and $z_0$, in fact it is of the form $r P_{n-1}(r,z_0)$ where $P_{n-1}(r,z_0)$ is a polynomial of degree $n-1$ in the variables $r$ and $z_0$. In particular, equation $\psi_n(r,z_0)=0$ defines an algebraic real curve of degree $n$ for each every solution $(\tilde{r}_0(\tilde{z}_0),\tilde{z}_0)$, satisfying $(\partial \psi_n/\partial r)(\tilde{r}_0(\tilde{z}_0),\tilde{z}_0)\neq0$, corresponds to a closed orbit in $\mathbb{R}^3$. That follows from Fundamental Lemma. Doing the procedure for every $z_0\in S$ we get $\mathfrak{M}$, so statements $(i)$ and $(ii)$ are proved.

The second part of the proof is done by induction on the degree $n$ of the perturbation. For $n=1$ we get that function $\psi_1(r,z_0)$ writes
$$
\psi_1(r,z_0)=r\left[\pi(a_{000}^{+}+a_{000}^{-})\right].
$$
Therefore the perturbation $g_3^+=g_1^-=g_2^-=g_3^-=0$ joint with
$$
\begin{array}{l}
g_1^+(x,y,z)=\dfrac{C_{00}^1}{\pi}x;\vspace{0.2cm}\\
g_2^+(x,y,z)=\dfrac{C_{00}^1}{\pi}y;
\end{array}
$$
realizes the polynomial $\bar{\psi}_1(r,z_0)= r\left[C_{00}^1\right]$. Assume that for any polynomial \[\widetilde{\psi}_n(r,z_0)=r \left[\displaystyle\sum_{0\leq i+j\leq n-2}\widetilde{C}_{ij}^nr^iz_0^j\right],\] 
of degree $n-1$, there exists a perturbation $\widetilde{g}^\pm=(\widetilde{g}_1^\pm,\widetilde{g}_2^\pm,\widetilde{g}_3^\pm)$, also of degree $n-1$, that realizes the polynomial $\widetilde{\psi}_n$. Let be given a polynomial of degree $n$ of the form \[\bar{\psi}_n(r,z_0)=r\left[\displaystyle\sum_{0\leq i+j\leq n-1}C_{ij}^nr^iz_0^j\right],\]
with arbitrary coefficients $C_{ij}^n$.
Consider the perturbation of degree $n$ given by
\begin{equation}\label{perturb}
\begin{array}{l}
g_1^+(x,y,z)= \widetilde{g}_1^+(x,y,z)+ x\sum_{i=0}^{n-1} a_{i}y^iz^{n-1-i},\\
g_2^+(x,y,z)= \widetilde{g}_2^+(x,y,z)+ y\sum_{i=0}^{n-1} a_{i}y^iz^{n-1-i}, \\
g_3^+(x,y,z)= \widetilde{g}_3^+(x,y,z), \\ 
g_1^-(x,y,z)= \widetilde{g}_1^-(x,y,z),  \\
g_2^-(x,y,z)= \widetilde{g}_2^-(x,y,z), \\ 
g_3^-(x,y,z)= \widetilde{g}_3^-(x,y,z).
\end{array}
\end{equation}
We will show how we can choose the coefficients $a_i$, of the perturbation, in terms of the coefficients $C^n_{ij}$ and $\widetilde{C}^n_{ij}$ of the polynomials $\bar{\psi}_n$ and $\widetilde{\psi}_n$, so that $g$ realizes the polynomial $\bar{\psi}_n$. From \eqref{perturb} it is clear that
\[ \begin{array}{ll}
\bar{\psi}_n(r,z_0) & = \displaystyle\widetilde{\psi}_n(r,z_0) + r\left[\int_0^\pi \sum_{i=0}^{n-1} a_ir^i(\sin\theta)^{i}(z_0+I_h(\theta))^{n-1-i}d\theta\right] \vspace{.3cm}\\
& \displaystyle=\widetilde{\psi}_n(r,z_0) + r\left[\sum_{i=0}^{n-1} \sum_{j=0}^{n-1-i} a_i \delta_{ij} r^iz_0^j\right],
\end{array}\]
where
$\delta_{ij} = \frac{(n-1-i)!}{j!(n-1-i-j)!} \int_0^\pi (\sin\theta)^{i}(I_h(\theta))^{n-1-i-j}d\theta.$
Here is important to mention that if $i+j=n-1$ then $\delta_{ij}\neq0$. So, we have the relation 
\[\left\{\begin{array}{ll}
C^n_{ij}=\widetilde{C}^n_{ij}+a_i\delta_{ij}, &\mbox{ for }0\leq i+j\leq n-2, \\
C^n_{ij}=a_i\delta_{ij}, &\mbox{ for }i+j= n-1.
\end{array}\right.\]
The induction works in the following way. Given the polynomial $\bar{\psi}_n$, first we take the coefficients $C^n_{ij}$ such that $i+j=n-1$, and we have the coefficients of the perturbation $a_i=C^n_{ij}/\delta_{ij}$. Now we consider the coefficients $\widetilde{C}^n_{ij}=C^n_{ij}-a_i\delta_{ij},$  for $0\leq i+j\leq n-2$, and by hypothesis of induction we obtain the perturbation $\widetilde{g}^\pm$. Finally, using \eqref{perturb}, we have the perturbation $g^\pm$ that realizes the polynomial $\bar{\psi}_n$. This concludes the proof of Theorem \ref{teo_generalv}.

\subsection{Proof of Theorem \ref{teo1}}

The first part of Theorem \ref{teo1} is a straightforward calculation using \eqref{function_F} and taking into account  that function $G$ in \eqref{final_system2} is piecewise linear in the variables $r$ and $z$. The second part of Theorem \ref{teo1} is an immediate consequence of statement $(ii)$ of Theorem \ref{teo_generalv}.

\subsection{Proof of Theorem \ref{teo2}}
The proof of Theorem \ref{teo2} is completely analogous to the proof of Theorem \ref{teo1} just observing that in this case the function $G$ in \eqref{final_system2} is piecewise quadratic in the variables $r$ and $z$.

\section{Some examples}\label{Examples}

In this section we present two examples, for Theorems \ref{teo1} and \ref{teo2}, in order to clarify the approach used throughout the paper.

\subsection{Example of Theorem \ref{teo1}}

Consider system \eqref{system_noncontinuous} putting $h(x,y)\equiv0$ and the radial perturbation \eqref{perturbations} with $a_{010}^{+}=1/2$, $a_{001}^{-}=1/\pi$ and all other coefficients vanishing. In this case system \eqref{system_noncontinuous} for $y\geq0$ writes
$$
\begin{array}{rcl}
\dot{x}&=&-y+\varepsilon x\left(\dfrac{y}{2}\right),\vspace{0.2cm}\\
\dot{y}&=&x-\varepsilon y\left(\dfrac{y}{2}\right),\vspace{0.2cm}\\
\dot{z}&=&0,
\end{array}
$$
and for $y\leq0$ it writes
$$
\begin{array}{rcl}
\dot{x}&=&-y+\varepsilon x\left(-\dfrac{z}{\pi}\right),\vspace{0.2cm}\\
\dot{y}&=&x+\varepsilon y\left(-\dfrac{z}{\pi}\right),\vspace{0.2cm}\\
\dot{z}&=&0,
\end{array}
$$
which is a piecewise quadratic system. Notice that since $\dot{z}=0$ every level $z=k$ with $k\in\mathbb{R}$ is invariant. We remark, as commented before, that every point on $\{y=0\}$ outside the $z-$axis is formed by crossing points.

\begin{figure}[t]
	\begin{center}
		\begin{overpic}[height=6cm]{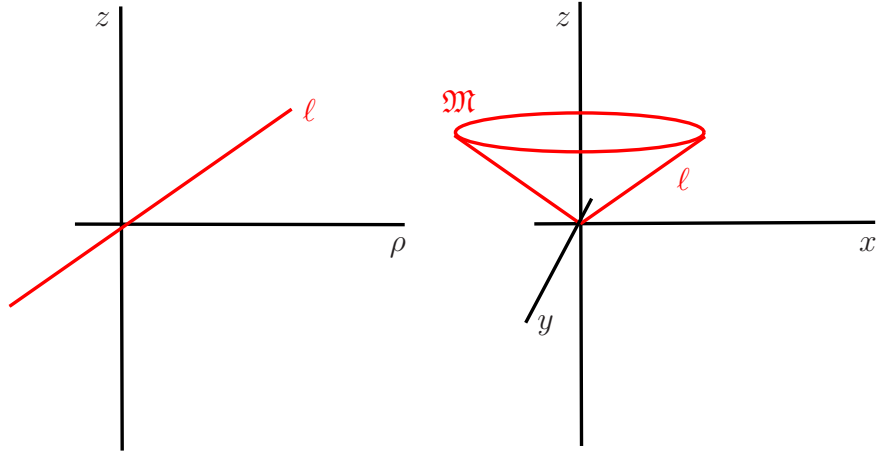} 
			\put(34,38){\textcolor{red}{$\ell$}}
			\put(77,30){\textcolor{red}{$\ell$}}
			\put(44,23){$\rho$}
			\put(10,49){$z$}
			\put(63,49){$z$}
			\put(98,23){$x$}
			\put(61,14){$y$}
			\put(50,39){\textcolor{red}{$\mathfrak{M}$}}
		\end{overpic}
		\caption{Example of Theorem \ref{teo1}.}\label{fig1}
	\end{center}
\end{figure}

The function $\psi_1(r,z_0)$ in this case is given by $\psi_1(r,z_0)=r(r-z_0)$, so $L_1$ is the straight line $r=z_0$ and $\ell=L_1\cap\{\theta=0\}\cap\{r>0\}$.
Now, since $C_{10}^1=1\neq0$ and $\ell$ is a half straight line (therefore non-empty) then $\widetilde{\mathfrak{M}}=\mathbb{S}^1\times\ell$ is a cone obtained by the revolution of $\ell$ around the $z-$axis, see Figure \ref{fig1}.

\subsection{Example of Theorem \ref{teo2}}

Consider system \eqref{system_noncontinuous} with $h(x,y)=x/\sqrt{x^2+y^2}$ which in cylindrical coordinates writes $h(r,\theta)=\cos\theta$. Clearly the function $h$ satisfies hypotheses ${\bf H1}$ and  ${\bf H2}$. Consider also perturbations \eqref{perturbations} with $a_{200}^+=2/\pi$, $a_{002}^+=1/(2\pi)$, $a_{010}^+=-3$, $a_{002}^-=1/(2\pi)$, $a_{000}^-=8/\pi$ and all other coefficients vanishing. Then system \eqref{system_noncontinuous} for $y\geq0$ writes
$$
\begin{array}{rcl}
\dot{x}&=&-y+\varepsilon x\left(\dfrac{2}{\pi}x^2+\dfrac{1}{2\pi}z^2-3y\right),\vspace{0.2cm}\\
\dot{y}&=&x+\varepsilon y\left(\dfrac{2}{\pi}x^2+\dfrac{1}{2\pi}z^2-3y\right),\vspace{0.2cm}\\
\dot{z}&=&\dfrac{x}{\sqrt{x^2+y^2}},
\end{array}
$$
and for $y\leq0$ it writes
$$
\begin{array}{rcl}
\dot{x}&=&-y-\varepsilon x\left(\dfrac{1}{2\pi}z^2+\dfrac{8}{\pi}\right),\vspace{0.2cm}\\
\dot{y}&=&x-\varepsilon y\left(\dfrac{1}{2\pi}z^2+\dfrac{8}{\pi}\right),\vspace{0.2cm}\\
\dot{z}&=&\dfrac{x}{\sqrt{x^2+y^2}}.
\end{array}
$$

In order to verify the other hypotheses of Theorem \ref{teo2} we note that in this particular case we get $I_h(\theta)=\sin\theta$. Due to the previous choice of parameters for the radial perturbation, we only have to calculate two integral involving $h$, namely, $c_{10}^0=2$ and $c_{01}^0=-2$. The other constants do not play any role because they are multiplied by zero. Thus function $\psi_2$ given in \eqref{poly} is $\psi_2(r,z_0)=r(r^2+z_0^2-6r+8)=r((r-3)^2+z_0^2-1)$.

For every $z_0$ satisfying $-3\leq z_0 \leq 4$, the zeros of $\psi_2$ belong to the circumference on $M$ with $\theta=0$ which is centered at the point $(\widetilde{r},\widetilde{z})=(3,0)$ having radius $\widetilde{\rho}=1$. Evaluating such zeros in $\partial\psi_2/\partial r$ we obtain $\mp2\sqrt{1-z^2}$ which is nonzero if $z\neq\pm1$. Therefore $S$ is the open interval $S=(-1,1)$, since $\psi_2$ has no zero if $|z_0|>1$.

Finally, since $\psi_2^{-1}(0)$ is a circumference (excluding north and south poles) on $[1,\rho_1]\times S$, call $\widetilde{\mathbb{S}}$, we distinguish the following cases:
\begin{itemize}
	\item If $\rho_1<2$, then $\psi_2^{-1}(0)=\emptyset$;
	\item If $2<\rho_1<4$, then $\psi_2^{-1}(0)$ is a subset of $\widetilde{\mathbb{S}}$ and $\widetilde{\mathfrak{M}}=\mathbb{S}^1\times\psi_2^{-1}(0)$ is a set contained on the tori $\widetilde{\mathfrak{T}}=\mathbb{S}^1\times\widetilde{\mathbb{S}}$, which is a tori where the parallels $\mathbb{S}^1\times\{(0,3,1)\}$ and $\mathbb{S}^1\times\{(0,3,-1)\}$ have been excluded.
	\item If $\rho_1>4$ then $\psi_2^{-1}(0)=\widetilde{\mathbb{S}}$ and $\widetilde{\mathfrak{M}}=\widetilde{\mathfrak{T}}$.
\end{itemize}

\begin{figure}
	\begin{center}
		\begin{overpic}[height=6cm]{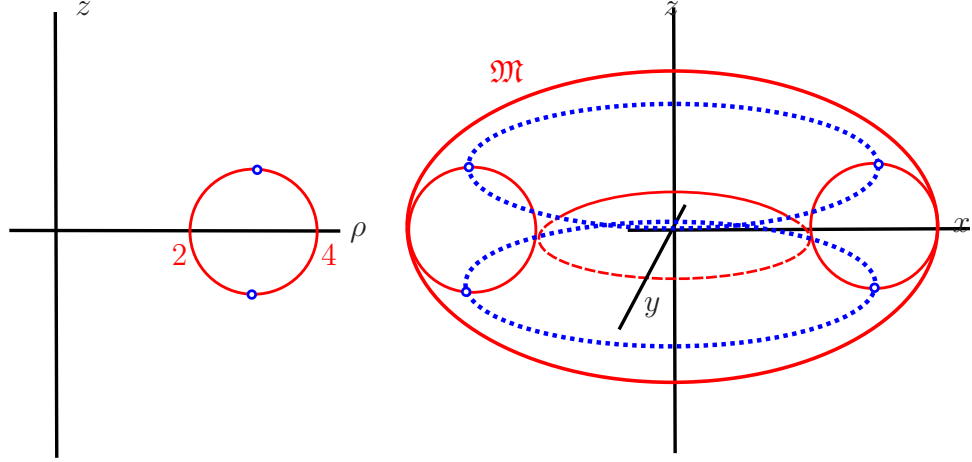} 
			\put(17,20){\textcolor{red}{$2$}}
			\put(32.5,20){\textcolor{red}{$4$}}
			\put(35.5,23){$\rho$}
			\put(7,46){$z$}
			\put(50,39){\textcolor{red}{$\mathfrak{M}$}}
			\put(68,46){$z$}
			\put(98,23.5){$x$}
			\put(66,15.3){$y$}
		\end{overpic}
		\caption{Example of Theorem \ref{teo2}.}\label{fig2}
	\end{center}
\end{figure}

Now if $\rho_1=2$ or $\rho_1=4$, then $\psi_2^{-1}(0)=\{p\}$ and therefore the method does not apply. See Figure \ref{fig2}.

\section{Averaging theory of first order for piecewise smooth
	differential systems}\label{section_averaging}

In this section we present the first-order averaging theory for piecewise smooth differential systems. This method was introduced in \cite{LliNovTei2015c} and is summarized in what follows.
\begin{theorem}\label{averaging}
	We consider the following piecewise smooth differential system
	\begin{equation}\label{MRs1}
	x'(t)=\e G(t,x)+\e^2R(t,x,\e),
	\end{equation}
	with
	\[
	\begin{aligned}
	&G(t,x)=G^+(t,x)+\sgn(\eta(t,x))G^-(t,x),\\
	&R(t,x,\e)=R^+(t,x,\e)+\sgn(\eta(t,x))R^-(t,x,\e),
	\end{aligned}
	\]
	where $G^+,G^-:\R\times D\rightarrow\R^n$, $R^+,R^-:\R\times
	D\times(-\e_0,\e_0)\rightarrow\R^n$ and $\eta:\R\times D\rightarrow \R$
	are continuous functions, $T$--periodic in the variable $t$ and $D$
	is an open subset of $\R^n$. We also suppose that $\eta$ is a $C^1$
	function having $0$ as a regular value. 
	
	Define the averaged function $\psi:D\rightarrow\R^n$ as
	\begin{equation*}\label{MRf1}
	\psi(x)=\int_0^T G(t,x) dt.
	\end{equation*}
	We assume the following three conditions.
	\begin{itemize}
		\item[(i)] $G^+,\,G^-,\,R^+,\,R^-$ and $\eta$ are locally Lipschitz
		with respect to $x$;
		
		\item[(ii)] there exists an open bounded subset $C \subset D$ such that, for $|\e|>0$ sufficiently small, every orbit starting in $\overline{C}$ reaches the set of piecewise smooth only at its crossing regions.
		
		\item[(iii)] for $a\in C$ with $\psi(a)=0$, there exist a neighborhood
		$U \subset C$ of $a$ such that $\psi(z)\neq 0$ for all $z\in\overline{U}
		\backslash\{a\}$ and $d_B(\psi,U,a)\neq 0$, (i.e. the Brouwer degree of
		$\psi$ at $a$ is not zero).
	\end{itemize}
	
	Then, for $|\e|>0$ sufficiently small, there exists a $T$--periodic
	solution $x(\cdot,\e)$ of system \eqref{MRs1} such that $x(0,\e)\to
	a$ as $\e\to 0$.
	
\end{theorem}

\begin{remark}\label{rem1}
In  \cite{LliNovTei2015c} is proved that if $\partial_t\eta(t,x)\neq 0$ for each $(t,x)\in\Sigma=\eta^{-1}(0)$ then hypothesis (ii) in Theorem \ref{averaging} holds.
\end{remark}

\begin{remark}
We observe that if function $\psi(z)$ is of class $C^1$, $\psi(a)=0$ and the Jacobian $J\psi(a)$ is not zero, then $d_B(\psi,V,a)\neq 0$. For more information concerning the Brouwer degree see \cite{Bro1983}.
\end{remark}

\section*{Acknowledgments}

The first author is partially supported by the Brazilian Capes grant 88881.068462/2014-01 and Brazilian FAPESP grant 2019/10269-3. The second author is partially supported by Pronex/FAPEG/CNPq  Proc. 2012 10 26 7000 803 and Proc. 2017 10 26 7000 508, Capes grant 88881.068462/2014-01 and Universal/CNPq grant 420858/2016-4. The third author is partially supported by Universal/CNPq grant 434599/ 2018-2 and Fapesp-Brazil 2018/13481-0.


\end{document}